\documentclass[a4paper,12pt,reqno]{amsart}
\usepackage[a4paper,margin=1.9cm,top=2.5cm,bottom=2.5cm,centering,vcentering]{geometry}
\usepackage{amssymb,mathtools,cite,enumerate,color,eqnarray,hyperref}
\definecolor{blue}{RGB}{0, 0, 200}
\definecolor{pink}{RGB}{252, 0, 50}

\hypersetup{
           colorlinks=true,
           linkcolor=blue,
           urlcolor=red,
           citecolor=pink,
          }

\theoremstyle{plain}
\newtheorem{theorem}{Theorem}[section]

\newtheorem{lemma}[theorem]{Lemma}

\theoremstyle{definition}

\newtheorem{defn}[theorem]{Definition}

\numberwithin{equation}{section}

\numberwithin{equation}{section}

\begin{document}

\title{Arithmetic of 2-regular partitions with distinct odd parts}

\author[H. Nath]{Hemjyoti Nath}
\address{Department of mathematical sciences, Tezpur University, Napaam, Tezpur, Assam 784028, India}
\email{hemjyotinath40@gmail.com}
\keywords{Integer partitions, Ramanujan-type congruences, Modular forms, Hecke eigenforms.}

\subjclass[2020]{05A17, 11P83, 11F11, 11F20.}

\date{\today.}

\begin{abstract}
Let $pod_2(n)$ denote the number of $2$-regular partitions of $n$ with distinct odd parts (even parts are unrestricted). In this article, we obtain congruences for $pod_2(n)$ mod $2$ and mod $8$ using some generating function manipulations and the theory of Hecke eigenform. 
\end{abstract}

\maketitle

\section{Introduction}
Let $p(n)$ be the number of partitions of a positive integer $n$ : the number of nonincreasing sequence of positive integers whose sum is $n$. By convention, we agree that $p(0):=1$. Partitions play an important role in number theory, combinatorics, representation theory, and mathematical physics. Euler showed that the generating function for $p(n)$ is
\begin{equation*}
    \sum_{n=0}^{\infty}p(n)q^n = \frac{1}{f_1},
\end{equation*}
where
\begin{equation*}
    f_j := \prod_{n=1}^{\infty}(1-q^{jn}), \quad |q|<1, \quad j \geq 1.
\end{equation*}
Among the most famous results on $p(n)$ are the Ramanujan congruences. For all nonnegative integers $n$, these are 
\begin{align*}
    p(5n+4)&\equiv 0 \pmod5, \\p(7n+5)&\equiv 0 \pmod7, \intertext{and} p(11n+7)&\equiv 0 \pmod{11}.
\end{align*}
\\
The generating function of $pod_2(n)$ is given by
\begin{equation}\label{e9}
    \sum_{n=0}^{\infty}pod_2(n)q^n = \frac{\psi(-q^2)}{\psi(-q)}
\end{equation}

Recently, the arithmetic properties of $pod_{\ell}(n)$ have been widely studied. See, for example, \cite{gir,hem,dre,vee,nai,JJ}. The purpose of this paper is to establish some Ramanujan-type congruences for $pod_2(n)$. In particular, we prove infinite families of congruences modulo $2$ and modulo $8$ for $pod_2(n)$. The following are our main results.

\begin{theorem}\label{t1}
For any positive integer $n$, we have
\begin{enumerate}[(i)]
    \item \begin{equation}\label{e10}
        pod_2(3n+2)\equiv 0 \pmod{2},
    \end{equation}
    \item \begin{equation}\label{e11}
        \sum_{n=0}^{\infty}pod_2(3n)\equiv f_1 \pmod{2},
    \end{equation}
    \item \begin{equation}\label{e12}
        \sum_{n=0}^{\infty}pod_2(9n+1)\equiv f_1^3 \pmod{2}.
    \end{equation}
\end{enumerate}    
\end{theorem}

Next, we present a complete characterization for the parity of $pod_2(3n+1)$.
\begin{theorem}\label{t5}
    Let $T_k$ denote the $k$-th triangular number, then for any positive integer $n$, we have
    \begin{equation}\label{t5.1}
        pod_2(3n+1) \equiv \begin{cases}
            1 \pmod{2}, \quad \text{if} \quad n=3T_k,\\
            0 \pmod{2}, \quad \text{otherwise}.
        \end{cases}
    \end{equation}
where $T_k = k(k+1)/2$.
\end{theorem}

In the next two theorems, we present some infinite families of congruences modulo $2$ and modulo $8$ respectively.

\begin{theorem}\label{t4}
Let $p$ be odd primes. Then the following statements hold.
\begin{enumerate}[(i)]
    \item Suppose that $s$ is an integer satisfying $1\leq s \leq 8p$, $s \equiv 1 \pmod{8}$ and $\left( \frac{s}{p} \right) = -1$, where $\left( \frac{s}{p} \right)$ is the Legendre symbol. Then,

\begin{equation}\label{e12.0.1}
    pod_2\left( pn+\frac{s-1}{8} \right) \equiv 0 \pmod{2}.
\end{equation}
If $\tau(p) \equiv 0 \pmod{2}$, then, for all $n \geq 0, k\geq 1$,
\begin{equation}\label{e12.0.2}
    pod_2\left( p^{2k+1}n+\frac{sp^{2k}-1}{8} \right) \equiv 0 \pmod{2}.
\end{equation}
    \item Suppose that $r$ is an integer  such that $1\leq r \leq 8p$, $rp \equiv 1 \pmod{8}$ and $(r,p)=1$. If $\tau(p) \equiv 0 \pmod{2}$, then, for all $n \geq 0, k\geq 1$,
    \begin{equation}\label{e12.0.3}
    pod_2\left( p^{2k+2}n + \frac{rp^{2k+1}-1}{8} \right) \equiv 0 \pmod{2}.
\end{equation}
\end{enumerate}
\end{theorem}

\begin{theorem}\label{t6}
Let $p$ be a prime such that $p \equiv 7 \pmod{8}$. Suppose that $r$ is an integer such that $1 \leq r < 8p$, $r \equiv 7 \pmod{8}$ and $(r,p)=1$. Then, for all $n \geq 0$ and $k \geq 0$, we have
\begin{equation}\label{et6}
    pod_2\left(p^{2k + 2}n + \frac{rp^{2k + 1} + 1}{8} \right) \equiv 0 \pmod{8}.
\end{equation}
\end{theorem}

We use Mathematica \cite{10} for our computations.\\

The paper is organised as follows: In Section \ref{pre}, we present some preliminaries required for our proofs. In Sections \ref{s1}--\ref{s4}, we  present the proofs of Theorems \ref{t1}--\ref{t6} respectively.

\section{Preliminaries}\label{pre}

The well known Ramanujan's general theta function $f(a,b)$ \cite[Equation 1.2.1]{6} is defined by
\begin{equation*}
    f(a,b)=\sum_{n=-\infty}^{\infty}a^{n(n+1)/2}b^{n(n-1)/2}, \quad |ab|<1.
\end{equation*}
Important special cases of $f(a,b)$ are the theta functions $\varphi(q)$, $\psi(q)$ and $f(-q)$, which satisfies the identities
\begin{equation*}
    \varphi(q) := f(q,q) = 1+2 \sum_{n=1}^{\infty}q^{n^2} = (-q;q^2)_{\infty}^2(q^2;q^2)_{\infty} = \frac{f_2^5}{f_1^2f_4^2},
\end{equation*}
\begin{equation*}
    \psi(q) := f(q,q^3) = \sum_{n=0}^{\infty}q^{n(n+1)/2} = \frac{(q^2;q^2)_{\infty}}{(q;q^2)_{\infty}} = \frac{f_2^2}{f_1},
\end{equation*}
and
\begin{equation*}
    f(-q):=f(-q,-q^2)=\sum_{n=-\infty}^{\infty}(-1)^{n}q^{n(3n+1)/2} = (q;q)_{\infty}=f_1.
\end{equation*}
From \textup{\cite[Equations 1.5.8, 1.5.9]{5}} we also note that
\begin{equation}\label{e9.0.0}
    \varphi(-q) = \frac{f_1^2}{f_2},
\end{equation}
and
\begin{equation}\label{e9.0}
    \psi(-q) = \frac{f_1f_4}{f_2}.
\end{equation}
In terms of $f(a,b)$, Jacobi's triple product identity \cite[Equation 1.3.11]{6} is given by
\begin{equation*}
    f(a,b) = (-a;ab)_{\infty}(-b;ab)_{\infty}(ab;ab)_{\infty}.
\end{equation*}

In order to prove theorems $\ref{t1}$, we require the following two lemmas.

\begin{lemma}\textup{\cite[Equation 14.3.3]{5}}\label{l1}
    We have
    \begin{equation}\label{e7}
        \frac{f_2^2}{f_1} = \frac{f_6f_9^2}{f_3f_{18}} + q\frac{f_{18}^2}{f_9}.
    \end{equation}
\end{lemma}

\begin{lemma}\textup{\cite[Equation 9]{11}}\label{l2}
    We have
    \begin{equation}\label{e8}
        \frac{f_2}{f_1^2} = \frac{f_6^4f_9^6}{f_3^8f_{18}^3} + 2q\frac{f_{6}^3f_9^3}{f_3^7}+4q^2\frac{f_6^2f_{18}^3}{f_3^6}.
    \end{equation}
\end{lemma}

To prove theorem \ref{t4}, we  will use the theory of Hecke eigenform. We describe some useful background material in the remainder of the section.

The Dedekind $\eta$- function is given by
\begin{equation*}
    \eta(z) := q^{1/24}(q;q)_{\infty},
\end{equation*}
where $q=e^{2\pi iz}$ and $z$ lies in the complex upper half plane $\mathbb{H}$. The well known $\Delta$-function is denoted by
\begin{equation*}
    \Delta(z) := \eta(z)^{24} = \sum_{n=1}^{\infty}\tau(n)q^n.
\end{equation*}

As usual, we denote $M_k(SL_2(\mathbb{Z}))$ (resp. $S_k(SL_2(\mathbb{Z}))$) is the complex vector space of weight $k$ holomorphic modular (resp. cusp) forms with respect to $SL_2(\mathbb{Z})$, where 
\begin{align*}
\text{SL}_2(\mathbb{Z}) & :=\left\{\begin{bmatrix}
a  &  b \\
c  &  d      
\end{bmatrix}: a, b, c, d \in \mathbb{Z}, ad-bc=1
\right\},
\end{align*}

\begin{defn}
    Suppose $f(z) = \sum_{n=0}^{\infty}a(n)q^n \in M_k(SL_2(\mathbb{Z})$, then for positive integer $m$, the action of Hecke operator $T_m$ on $f(z)$ is defined by
\begin{align*}
f(z)|T_{m,k} := \sum_{n=0}^{\infty} \left(\sum_{d\mid \gcd(n,m)}d^{k-1}a\left(\frac{nm}{d^2}\right)\right)q^n.
\end{align*}
In particular, if $m=p$ is prime, we have 
\begin{align}\label{hecke1}
f(z)|T_{p,k} := \sum_{n=0}^{\infty} \left(a(pn)+p^{k-1}a\left(\frac{n}{p}\right)\right)q^n.
\end{align}
We note that $a(n)=0$ unless $n$ is a nonnegative integer.
\end{defn}

\begin{defn}\label{hecke2}
A modular form $f(z)=\sum_{n=0}^{\infty}a(n)q^n \in M_k(SL_2(\mathbb{Z}))$ is called a Hecke eigenform if for every $m\geq2$ there exists a complex number $\lambda(m)$ for which 
\begin{align*}
f(z)|T_m = \lambda(m)f(z).
\end{align*}
\end{defn}

It is known that $\Delta(z) \in S_{12}(SL_2(\mathbb{Z}))$ and is an eigenform for all Hecke operators. In other words
\begin{equation*}
    \Delta(z) \mid T_{p,12} = \tau(p)\Delta(z)
\end{equation*}
for any prime $p$. The coefficients of $\Delta(z) = \sum_{n=1}^{\infty}\tau(n)q^n$ satisfies the following properties
\begin{equation}\label{tau1}
    \tau(mn) = \tau(m)\tau(n), \quad \text{if} \quad gcd(m,n)=1,
\end{equation}
\begin{equation}\label{tau2}
    \tau(p^{\ell}) = \tau(p)\tau(p^{\ell-1}) - p^{11}\tau(p^{\ell-2}).
\end{equation}

\section{Proof of theorem \ref{t1}}\label{s1}

In this section, we prove some congruences modulo $2$, using some generating function manipulations.\\

We have from $\eqref{e9}$ and $\eqref{e9.0}$ that
\begin{equation}\label{e1.0.0.0}
    \sum_{n=0}^{\infty}pod_2(n)q^n = \frac{f_2^2f_8}{f_1f_4^2}.
\end{equation}
Employing $\eqref{e7}$ and $\eqref{e8}$ in $\eqref{e1.0.0.0}$, we get
\begin{align}
    \sum_{n=0}^{\infty}pod_2(n)q^{n} & = \frac{f_6f_9^2f_{24}^4f_{36}^6}{f_3f_{18}f_{12}^8f_{72}^3} + 2q^4\frac{f_6f_9^2f_{24}^3f_{36}^3}{f_3f_{18}f_{12}^7} + 4q^8\frac{f_6f_9^2f_{24}^2f_{72}^3}{f_3f_{18}f_{12}^6} \nonumber \\
    & + q\frac{f_{18}^2f_{24}^4f_{36}^6}{f_9f_{12}^8f_{72}^3} + 2q^5\frac{f_{18}^2f_{24}^3f_{36}^3}{f_{9}f_{12}^7} + 4q^9\frac{f_{18}^2f_{24}^2f_{72}^3}{f_9f_{12}^6}\label{e13}.
\end{align}

Extracting the terms involving $q^{3n+2}$ from both sides of $\eqref{e13}$, dividing both sides by $q^2$ and then replacing $q^3$ by $q$, yields

\begin{equation}\label{e14}
    \sum_{n=0}^{\infty}pod_2(3n+2)q^{n} = 2q\frac{f_6^2f_8^3f_{12}^3}{f_3f_4^7} + 4q^2\frac{f_2f_3^2f_8^2f_{24}^3}{f_1f_6f_4^6}.
\end{equation}

Therefore, it follows from $\eqref{e14}$ that 
\begin{equation*}
   pod_2(3n+2) \equiv 0 \pmod {2}.
\end{equation*}
Thus, congruence $\eqref{e10}$ is obtained.\\

Next, extracting the terms involving $q^{3n}$ from both sides of $\eqref{e13}$ and then replacing $q^3$ by $q$, yields

\begin{equation}\label{e15}
    \sum_{n=0}^{\infty}pod_2(3n)q^{n} = \frac{f_2f_3^2f_8^4f_{12}^6}{f_1f_6f_4^8f_{24}^3}+4q^3\frac{f_6^2f_8^2f_{24}^3}{f_3f_4^6}.
\end{equation}

But, by the binomial theorem, $f_t^{2m} \equiv f_{2t}^m \pmod{2}$, for all positive integers $t$ and $m$. 

Therefore, it follows from $\eqref{e15}$ that 
\begin{equation*}
    \sum_{n=0}^{\infty}pod_2(3n)q^{n} \equiv f_1 \pmod {2}.
\end{equation*}

Thus, congruence $\eqref{e11}$ is obtained.\\

Next, extracting the terms involving $q^{3n+1}$ from both sides of $\eqref{e13}$, dividing both sides by $q$ and then replacing $q^3$ by $q$, yields

\begin{equation}\label{e16}
    \sum_{n=0}^{\infty}pod_2(3n+1)q^{n} = \frac{f_{6}^2f_{8}^4f_{12}^6}{f_3f_{4}^8f_{24}^3}+2q\frac{f_2f_3^2f_{8}^3f_{12}^3}{f_1f_{4}^7f_{6}}.
\end{equation}

Therefore, it follows from $\eqref{e16}$ that 
\begin{equation}\label{e17}
    \sum_{n=0}^{\infty}pod_2(3n+1)q^{n} \equiv f_3^3 \pmod {2}.
\end{equation}

Extracting the terms involving $q^{3n}$ from both sides of $\eqref{e17}$ and then replacing $q^3$ by $q$, yields

\begin{equation}\label{e18}
    \sum_{n=0}^{\infty}pod_2(9n+1)q^{n} \equiv f_1^3 \pmod {2}.
\end{equation}

Thus, congruence $\eqref{e12}$ is obtained.\\

\section{Proof of theorem \ref{t5}}\label{s2}
We have from $\eqref{e17}$ that
\begin{equation*}
    \sum_{n=0}^{\infty}pod_2(3n+1)q^{n} = f_3^3 \pmod {2}.
\end{equation*}
By the Jacobi's identity \textup{\cite[Equation 1.3.24]{6}}, we have for $|q|<1$,
\begin{equation}\label{e63}
   f_1^3 = \sum_{n=0}^{\infty}(-1)^n(2n+1)q^{n(n+1)/2}.
\end{equation}
Therefore,
\begin{equation}\label{e63.1}
   f_1^3 \equiv \sum_{n=0}^{\infty}q^{n(n+1)/2} \pmod{2}.
\end{equation}
Changing $q \to q^3$ in the above equation, we arrive at $\eqref{t5.1}$.

\section{Proof of theorem \ref{t4}}\label{s3}
In this section, we prove the infinite family of congruences modulo $2$ using the theory of Hecke eigenforms.\\

Magnifying $q \to q^8$ and multiplying $q$ on both sides of $\eqref{e1.0.0.0}$, we have
\begin{equation}\label{e4}
    \sum_{n=0}^{\infty}pod_2(n)q^{8n+1}=q\frac{f_{16}^2f_{64}}{f_{8}f_{32}^2} \equiv \Delta(z) = \sum_{n=1}^{\infty}\tau(n)q^n \pmod{2}.
\end{equation}
Hence, from $\eqref{e63.1}$, we have
\begin{equation*}
    \Delta(z)=qf_1^{24} \equiv qf_8^3 \equiv \sum_{n=0}^{\infty}q^{(2n+1)^2} \pmod{2}.
\end{equation*}
Therefore, we have
\begin{equation}\label{e5}
    \sum_{n=0}^{\infty}pod_2(n)q^{8n+1}\equiv \sum_{n=0}^{\infty}q^{(2n+1)^2} \pmod{2}.
\end{equation}
If $s\equiv1 \pmod{8}$ and $\left(\frac{s}{p}\right) = -1$, then, for any $n\geq 0$, $8np+s$ cannot be an odd square. This implies that the coefficients of $q^{8np+s}$ in the left-hand side of $\eqref{e5}$ must be even. It follows that
\begin{equation}\label{e6}
    pod_2\left( pn + \frac{s-1}{8} \right) \equiv 0 \pmod{2}.
\end{equation}
This proves $\eqref{e12.0.1}$.\\

Since $\tau(p) \equiv 0 \pmod{2}$ and $\Delta(n)$ is a Hecke eigenform, we have
\begin{equation*}
    \Delta(z)\mid T_{p,12} = \tau(p)\Delta(z) \equiv 0 \pmod{2}.
\end{equation*}
By $\eqref{hecke1}$ and $\eqref{e4}$, we get
\begin{equation*}
    \sum_{n=0}^{\infty}pod_{2}(n) q^{8n+1} \mid T_{p,12} \equiv \sum_{n = n_0}^{\infty}\left( pod_2\left( \frac{pn-1}{8} \right) + pod_2\left( \frac{n/p-1}{8}\right) \right)q^n \equiv 0 \pmod{2}.
\end{equation*}
If we write $m=\frac{n/p-1}{8}\in N $, then $\frac{pn-1}{8} = p^2m + \frac{p^2-1}{8}$. So, we have
\begin{equation*}
    pod_2\left( p^2m + \frac{p^2-1}{8} \right) + pod_2(m) \equiv 0 \pmod{2}.
\end{equation*}
That is, 
\begin{equation*}
    pod_2\left( p^2m + \frac{p^2-1}{8} \right) \equiv pod_2(m) \pmod{2}.
\end{equation*}
By induction, for $k\geq 1$ we find that
\begin{equation*}
    pod_2\left( p^{2k}m + \frac{p^{2k}-1}{8} \right) \equiv pod_2(m) \pmod{2}.
\end{equation*}
Considering $\eqref{e6}$, we get
\begin{equation*}
    pod_2\left( p^{2k+1}n+\frac{sp^{2k}-1}{8} \right) \equiv pod_2\left( pn + \frac{s-1}{8} \right) \equiv 0 \pmod{2}.
\end{equation*}
This proves $\eqref{e12.0.2}$.\\

Since $\tau(p) \equiv 0 \pmod{2}$, formula $\eqref{tau2}$ gives
\begin{equation*}
    \tau(p^{2k+1}) \equiv 0 \pmod{2}.
\end{equation*}
Applying $\eqref{tau1}$, we have
\begin{equation*}
    \tau\left( p^{2k+1} \left( 8pn + r \right) \right) = \tau\left( p^{2k+1} \right)\tau\left( 8pn+r \right) \equiv 0 \pmod{2}.
\end{equation*}
It follows from $\eqref{e4}$ that
\begin{equation*}
    pod_2\left( p^{2k+2}n + \frac{rp^{2k+1}-1}{8} \right) \equiv 0 \pmod{2}.
\end{equation*}
This proves $\eqref{e12.0.3}$.\\

\section{Proof of theorem \ref{t6}}\label{s4}
In this section, we prove the infinite family of congruences modulo $8$.\\

We have
\begin{align*}
    \sum_{n=0}^{\infty}pod_2(n)q^n & = \frac{f_2^2f_8}{f_1f_4^2}\\
    & = \frac{\psi(q)}{\varphi(-q^4)}\\
    & = \frac{\psi(q)}{1+ (\varphi(-q^4)-1)}\\
    & = \psi(q)\left( 1-(\varphi(-q^4)-1) + (\varphi(-q^4)-1)^2- \ldots \right),
\end{align*}

then from modulo $8$, we have

\begin{align*}
    \sum_{n=0}^{\infty}pod_2(n)q^n & \equiv \psi(q)\left( 1-(\varphi(-q^4)-1) + (\varphi(-q^4)-1)^2 \right) \pmod{8}\\
    & \equiv \psi(q)\left( 1-2\sum_{n=1}^{\infty}(-1)^nq^{4n^2} + 4\left( \sum_{n=1}^{\infty}(-1)^nq^{4n^2} \right)^2 \right) \pmod{8}.
\end{align*}

Now
\begin{align}\label{e19}
    \left( \sum_{n=1}^{\infty}(-1)^nq^{n^2} \right)^2 & = \sum_{n=1}^{\infty}q^{2n^2} + 2\sum_{1 \leq n_1 < n_2}(-1)^{n_1+n_2}q^{n_1^2+n_2^2}\\
    & \equiv \sum_{n=1}^{\infty}q^{2n^2} \pmod{2},
\end{align}

so modulo $8$

\begin{align}\label{e20}
    \sum_{n=0}^{\infty}pod_2(n)q^n \equiv \psi(q)\left( 1-2\sum_{n=1}^{\infty}(-1)^nq^{4n^2} + 4\sum_{n=1}^{\infty}(-1)^nq^{8n^2} \right) \pmod{8}.
\end{align}

Now we have
\begin{align}\label{e21}
    \sum_{n=0}^{\infty}pod_2(n)q^{8n+1} \equiv \left( 1-2\sum_{n=1}^{\infty}(-1)^nq^{32n^2} + 4\sum_{n=1}^{\infty}(-1)^nq^{64n^2} \right) \times \left( \sum_{n=0}^{\infty}q^{(2n+1)^2} \right) \pmod{8}.
\end{align}

Note that if $8n+1 = p^{2\alpha+1}(8pm+r)$ where $p\equiv 7 \pmod{8}, (r,p)=1, rp \equiv 1 \pmod{8},$ then $8n+1$ cannot be represented by $x^2$, $x^2+y^2$ or $x^2+2y^2$ (see \cite{4}). So for these $8n+1$, from $\eqref{e20}$ and $\eqref{e21}$ it is easy to see that $pod_2(n) \equiv 0 \pmod{8}$. This completes the proof.


\begin{thebibliography}{99}

\bibitem{6}
B.C. Berndt, \textit{Number Theory in the Spirit of Ramanujan}, AMS, Providence, 2006.

\bibitem{hem}
B. Hemanthkumar, H. S. Sumanth Bharadwaj, and M. S. Mahadeva Naika, 
Arithmetic properties of 9-regular partitions with distinct odd parts.
\textit{Acta Mathematica Vietnamica},
\textbf{44}, (2019), 797–-811.

\bibitem{gir}
D. S. Gireesh, M. D. Hirschhorn, and M. S. Mahadeva Naika,                        
On 3-regular partitions with odd parts distinct.
\textit{Ramanujan J.}, 
\textbf{44}, (2017), 227–-236.

\bibitem{7}
H. Nath,
The pod function and its connection with other partition functions.
\textit{Rad HAZU, Matematičke znanosti},
(2024), (to appear).

\bibitem{8}
H. Nath,
Parity results of PEND partition. preprint \textbf{arXiv:2407.10428},
(2024).

\bibitem{9}
H. Nath,
New congruences for Partitions where the Even Parts are Distinct. 
\textit{Integers},
(2024), (to appear).

\bibitem{JJ}
Jing-Jun Yu, 
Some congruences for regular partitions with distinct odd parts.
\textit{Integers},
\textbf{23}, (2023), \#A77.

\bibitem{5}
M.D. Hirschhorn, \textit{The power of \lowercase{$q$}, a personal journey}, 
Developments in Mathematics, 49.
Springer, 2017.

\bibitem{4}
M.D. Hirschhorn,
Partial fractions and four classical theorems of number theory. 
\textit{Amer. Math. Monthly},
\textbf{107}, (2000), no. 3, 260--264.

\bibitem{nai}
N. Saikia, 
Infnite families of congruences for 3-regular partitions with distinct odd parts. 
\textit{Commun. Math. Stat.},
\textbf{8}, (2020), no. 4, 443–-451.

\bibitem{dre}
R. Drema and N. Saikia, 
Arithmetic properties for $\ell$-regular partition functions with distinct even parts.
\textit{Bol. Soc. Mat. Mex.},
\textbf{28}, (2022), 10–-20.

\bibitem{11}
R. Guadalupe,
Some congruences for $3$-core cubic bipartitions.
preprint \textbf{arXiv:2311.17674},
(2023).

\bibitem{vee}
V. S. Veena and S. N. Fathima, 
Arithmetic properties of 3-regular partitions with distinct odd parts.
\textit{Abh. Math. Semin. Univ. Hambg.},
\textbf{91}, (2021), 69–-80.

\bibitem{10}
Wolfram~Research{,} Inc.
\newblock Mathematica, {V}ersion 10.0.
\newblock Champaign, IL, 2014. 

\end{thebibliography}
\end{document}